\documentclass[12pt]{article}
\usepackage{amsmath,amssymb,amsthm}
\usepackage[xdvi]{graphicx}
\usepackage[T1]{fontenc}
\usepackage[applemac]{inputenc}
\newcommand{\<}{\langle}
\newcommand{\argmax}{\operatorname*{argmax}}
\renewcommand{\>}{\rangle}

\newcommand{\mma}{\texttt{mathematica}}
\newcommand{\moy}{\operatorname{moy}}
\newcommand{\str}{\operatorname{strat}}
\renewcommand{\vec}{\mathbf}
\title{Fate stochastic management and policy benchmark in 421, a popular game}
\author{Pierre \textsc{Albarède}\\
\small b. A, r\'es. Valvert, 12, av. de la Fourane, 
13090 Aix-en-Provence, France\\
\small palbarede@yahoo.com}
\date{\today}
\begin{document}
\maketitle
\begin{abstract}
    Using game and probability theories, I study the French popular game 
    421, a perfect information stochastic stage game.  The problem is 
    to find strategies maximizing the probability of some expected 
    utility.  I only solve a player's round against providence, a 
    problem of fate stochastic management: beyond the backward 
    induction solution, bounded complexity motivates heuristic 
    policies.  For a unique goal utility, a simple optimal policy, 
    ratchet, is obtained.  Its result probabilities are compiled and 
    used, for arbitrary utilities, as the logic of goal identification 
    policies.  Various policies appear, close to human behavior, and 
    are exactly evaluated by solving the Kolmogorov equation.

\noindent k.  w.: stochastic management, Kolmogorov equation, bounded 
complexity, human behavior.

\noindent JEL  
C61, 
C63, 
C73. 
%
%
%
%
%
%
%
%
%
%
%
%
\end{abstract}
\newpage
\tableofcontents
\newpage
\section{Aim and interest of the study}

Following \cite{421MODE}, I look for strategies, maximizing the 
probability of win, or some expected utility, in the game 421 
combining chance and decision (see appendix \ref{reg}).

By the way, or, indeed, by \emph{serendipity}, I encounter the problem 
of fate \emph{stochastic management}: optimizing today's decisions, 
with respect to a future utility, and in spite of tomorrow's odds.  
Such issues as what are the optimal policies, in what circumstances, 
and how much they demand on intellectual resources, can be resolved 
mathematically, suggesting that management could be an exact science 
(as part of operations research).

A lottery is not a game, in the sense of game theory, but a stochastic 
process (a sequence of random variables).  Game theory treats 
classically multiple player decision games, the archetype of which is 
chess.  A game in which players' fates depend on both chance and their 
decisions, like 421 and backgammon \cite{Tesauro}, is a 
\emph{stochastic game} \cite{Shapley}.  Chance makes decision more 
complex.  For example, consider a variation on chess: the player 
proposes a list of $n$ moves, the actual one being determined by 
casting dice.  $n=1$ yields the standard pure decision game; 
$n=n_{1}$, where $n_{1}$ is the number of possible moves, yields a 
pure chance game; $n\approx n_{1}/2$ yields a game of chance and 
decision, more complex than the former and the latter.

Game theory primarily focuses on the existence proofs for optimal 
strategies.  However, ``usable techniques for obtaining practical 
answers'' also matter \cite[\S 1.1]{Isaacs}.  Indeed, little can be 
done from existence without construction: this is the old debate 
around Zermelo's axiom of choice.  Hence the interest of 
investigating, as in Church's thesis \cite{Delahaye}, calculability, 
the existence of an algorithmic solution.  But even calculability may 
not be sufficient for actual computation.  For example, consider again 
chess, a finite but very large game: the algorithmic solution provided 
by the Zermelo theorem \cite[ch.  6]{Berge} is of no practical use 
(until the final moves), as noticed by \cite[§11.4]{Dutta}, as it 
exceeds the capacity of any computer.  The study of finite games does 
not stop with Zermelo theorem, and this is because of complexity 
boundedness.  Algorithms shall be compared not only with respect to 
optimality (degree of completion of the task) but also complexity, 
using a bit of \emph{complexity theory} \cite{Delahaye}.

An algorithm is characterized by its optimality, size and computing 
time on a given computer, specialized by high-level functions and 
data.  The algorithm may be good or bad, short or long, fast or slow.  
The three qualities and quantities are not independent: the exchange 
of computing time against size is the principle of data compression, 
the exchange of computing time against optimality is the principle of 
heuristics.  When setting two quantities, the minimum of the third 
one, as a function of other implicit parameters, can be defined: the 
minimum size and the minimum computing time are respectively related 
with the Kolmogorov complexity and the Bennett logical depth 
\cite{Delahaye}.

For strategy-generating algorithms, or deciding algorithms, or 
\emph{policies}, optimality is the expected utility.  For a perfect 
information stage game, it is interesting, for possible extensions, to 
characterize the asymptotic behavior of the computing time, when the 
depth tends to infinity, i.  e. to know whether the algorithm is 
linear or polynomial, rather than exponential, as feared from the tree 
structure of the game.

The present study is thus an occasion to relate with each other, on a 
live case, various tools and concepts attached to games, processes, 
probabilities, control, programming, algorithmics and complexity, with 
applications in management and game practice.

\section{Backward induction optimal policy\label{backind}}

A stochastic game reduces to a pure decision game, by considering 
providence as a particular player \cite[ch.  4]{Guerrien}, whose mixed 
strategy, known a priori, results from usual statistical postulates 
(independence, stationarity) and cannot be 
optimized.\footnote{Probability theory began as the study of the 
providential strategy in chance games, at the time of Bayes or the 
Bernoullis.} 421, thus reduced, and with some precautions on the rules 
(appendix \ref{reg}), is a perfect information finite game and the 
Zermelo theorem applies.

I will solve only a sub-game, the player's round against providence 
(while other players stand still), a stochastic management problem, 
featuring a martingale problem and, for the first player, a stopping 
time problem \cite{Hunt}.  The analogy with Brownian motion provides 
statistical mechanics tools.

\subsection{Alea}

Let $D\in\mathbb{N}$ be the number of dice, normally 3, and 
$F\in\mathbb{N}^{*}$ the number of faces of every dice, normally 6.  
Dice are discernible\footnote{Discernibility is not an innocuous 
hypothesis, as shown by Gibbs' paradox \cite{Reif}.}, so that the 
probability space is the set of face sequences, or arrangements.  The 
class of arrangements corresponding to each other by a permutation is 
a combination, e.  g., ``n\'enette'', 221, is the subset of 
arrangements $\{(1,2,2),(2,1,2),(2,2,1)\}$, of redundancy three.  In 
421, which dice produced which face does not matter, because ranking 
depends only on combination.

I describe the die system as in statistical mechanics: each die is a 
particle, with only one \emph{phase} variable, \emph{face}.  The laws 
of mechanics are replaced by usual statistical hypothesis, abstracting 
chance from any specific random generator.  The system is described, 
in Lagrangian notation, by a face combination, or, in Eulerian 
notation, by the sequence $d_{f}$ of occupation numbers of every face 
$f=1\ldots F$, e.  g. the Lagrangian notation 421 translates into the 
Eulerian notation $(1,1,0,1,0,0)$ ($F=6$).

The interest of Eulerian notation lies in that the set of Eulerian 
combinations \[\vec{d}=(d_{f})_{f=1\ldots F}.\] is the partially 
ordered normed vector space $\mathbb{Z}^{F}$.  The canonic basis 
$\vec{e}_{g}=(\delta_{f,g})_{f}$ is aligned with ``brelans'', 
combinations with all faces of a kind.  I define the ball 
\[B(D\in\mathbb{N})\equiv\{\vec{d}\in\mathbb{Z}^{F},|\vec{d}|\le D\}\] 
and similarly (replacing $\le$ above by $=$ or $<$ ) the sphere 
$\partial B(D)$ and the open ball $\breve{B}(D)$.  The intersections 
with the positive cone are represented by $+$ exponents; the set of 
actual combinations is $B^{+}(D)$.  The norm of a combination is the 
sum of Eulerian component absolute values.  The norm of a positive 
combination is just its number of dice.  The canonic order $\le$, 
partial on $\mathbb{Z}^{F}$, differs from the hierarchic order 
$\preceq$ (\ref{ordhie}), total on $\partial B^{+}(D)$.

Distinct casts are independent and the probability of any face to be 
on top is $1/F$ (unloaded dice).  The arrangements of one combination 
are thus equiprobable, and the probability of a combination is just 
that of any of its arrangements, times the combination redundancy.  
For example, the probability of obtaining the combination 21 is 
$2/F^{2}$, while the probability of obtaining the combination 11 is 
$1/F^{2}$.  More generally, the probability of obtaining the 
combination $\vec{d}$, after one cast, is given by the multinomial 
law, with usual notations generalizing power and factorial to integer 
vectors:
\begin{equation}
    p(\vec{d})=\vec{p}^{\vec{d}}\frac{|\vec{d}|!}{\vec{d}!}, 
    \vec{p}=\frac{1}{F}(1\ldots 1)\in\mathbb{Q}^{F}, 
    \sum_{\vec{d}\in\partial B^{+}(D)}p(\vec{d})=1.
    \label{multinomial}
\end{equation}

\subsection{Fate}

For all $j\in\mathbb{N}$, let the state $\vec{d}_{j}$ be the 
combination, accumulated after $j$ casts, and the event 
$\vec{d}_{j+1/2}$ be the combination, obtained from the $j+1$-th cast.  
Fate is the infinite state and event alternate sequence
\begin{equation}
    \varphi\equiv(\vec{d}_{0},\vec{d}_{1/2},\vec{d}_{1},\vec{d}_{3/2}\ldots).
    \label{destin}
\end{equation}
The integer or half-integer index is used as a discrete time, integer 
time for states, half-integer for events.  The set of possible fates 
is described by the fate tree, where branching represents chance (from 
integer time to half-integer time) or decision (conversely).

The rules of 421 imply:
\begin{eqnarray}
    \vec{d}_{0}&\equiv&\vec{0},
    \label{d0=0}\\
    \forall j\in\mathbb{N},\vec{d}_{j+1/2}&\in&B^{+}(D_{j}\equiv D-|\vec{d}_{j}|), 
    \label{toutrejouer}\\
    0\le\vec{d}_{j+1}-\vec{d}_{j}&\le&\vec{d}_{j+1/2},
    \label{mortouvif}\\
    \exists(j\in\mathbb{N},j\le J),\vec{d}_{j}&\in&B^{+}(D),
    \label{stop1}
\end{eqnarray}
where $J\in\mathbb{N}$ is the maximum round duration, normally 3.  
$D_{j}$ is the number of live dice, which have not been accumulated 
after $j$ events and one state.  From (\ref{d0=0}, \ref{stop1}), 
$D_{0}=D,D_{J}=0$.

From (\ref{toutrejouer}, \ref{mortouvif}),
\begin{eqnarray}    
    D_{j+1}&\le&D_{j},
    \label{sort}\\  
    (\vec{d}_{j+1}-\vec{d}_{j}=\vec{d}_{j+1/2}) 
    &\Leftrightarrow&\vec{d}_{j+1}\in B^{+}(D).
    \label{tout}
\end{eqnarray}

The effective round duration $J_{1}$ is the minimum of $j$ in 
(\ref{stop1}).  The next players' effective round durations must equal 
the first player's.  Therefore, for all players,
\begin{eqnarray}
    \forall(j\in\mathbb{N},j<J_{1}),\vec{d}_{j}\in\breve{B}^{+}(D)&,&
    \vec{d}_{j+1/2}\neq 0,
    \label{stop2a}\\
    \vec{d}_{J_{1}}\in B^{+}(D)&,&
    \label{stop2b}\\
    \forall(j\in\mathbb{N},j>J_{1}),\vec{d}_{j-1/2}=0&,&
    \vec{d}_{j}=\vec{d}_{J_{1}}.
    \label{asymp}
\end{eqnarray}
(\ref{stop2a}, \ref{stop2b}) are used, firstly, after the first 
player's end of round, to determine $J_{1}$, subsequently, as 
additional rules for next players.  When $j$ increases, the state 
vector $\vec{d}_{j}$ moves in the positive ball, off the origin, 
towards its boundary where it gets stuck at $\vec{d}_{J_{1}}$, the 
round \emph{result}.  Fate is virtually continued by an infinite 
sequence, asymptotically alternating the result and the null event.

\subsection{Utility\label{uti}}

Following von Neumann and Morgenstern \cite[ch.  27]{Dutta}, a 
player's utility is a number, given by a causal function, i.  e. a 
function of history (past fate), compatible with the player's 
preferences, and such that the utility before a random event is just 
the expected utility, i.  e. the probability-weighted utility average, 
over possible outcomes.  Thus, expected utility is anti-causal, i.  e. 
prescribing utilities at some future time determines its expectation 
at all prior times.

One never knows when a game actually stops, as a it is often embedded 
in a larger game.  Tennis is a familiar example: a tennis ``game'' is 
a sub-game of a set, itself a sub-game of a match, tournament, ranking 
system\ldots this cascade does not even stop with a player's life, 
because of cooperation between individuals.  But, if we want to obtain 
any result, we must stop somewhere in the game cascade, and judge 
utility more or less empirically.  (Quite similarly, in mechanics or 
thermodynamics, the studied system is coupled with the rest of the 
world, by an often delicate boundary or cut-off condition.)

The study of 421 should stop at end of game, by setting players' 
utilities, for example, a binary utility: one for win, zero for loss, 
or incorporating economy, \`a la Bernoulli, the logarithm of earning 
divided by wealth \cite{Jaynes}.  However, I treat only the round.  At 
end of round, the Bernoulli formula does not make sense and utility is 
not given directly by the rules (in particular, the transfer function 
of table \ref{transf}).  By examining the rules, a few properties of 
utility are obtained; for example, at constant time, for a rational 
player, utility must be compatible with the hierarchic order 
(\ref{ordhie}), etc.
 
But I will not further characterize utility.  On the contrary, I will 
consider the round independently of the rest of the game, with 
arbitrary utilities, in order to treat the problem of fate stochastic 
management in a rather general way.

For all fate $\varphi$ (\ref{destin}), utility is judged at some time 
$J_{\varphi}$, either integer or half-integer in general (in the 
round, $J_{\varphi}\in\{1/2,3/2,5/2\}$), as a causal function:
\begin{equation}
    u(\vec{d}_{0},\vec{d}_{1/2}\ldots\vec{d}_{J_{\varphi}})\in\mathbb{Q}.
    \label{causale}
\end{equation}
The function $u$ has a variable number of arguments, formally, it is 
defined on $\bigcup_{j\in\mathbb{N}^{*}}B^{+}(D)^{j}$.  Utility is 
judged forever:
\begin{equation}
    u(\ldots\vec{d}_{J_{\varphi}},.)
    \equiv u(\ldots\vec{d}_{J_{\varphi}}).
    \label{udd=ud}
\end{equation}

The rules (\ref{toutrejouer}, \ref{mortouvif}, \ref{stop1}, 
\ref{stop2a}, \ref{stop2b}) are superseded by $-\infty$ utilities for 
rule breaking histories (excluding cheating).  In particular, the 
utilitarian version of (\ref{toutrejouer}) is
\begin{equation}
    \forall(j\in\mathbb{N},\vec{d}_{j+1/2}\notin B^{+}(D_{j})),
    u(\ldots\vec{d}_{j},\vec{d}_{j+1/2})=-\infty
    \label{mortouvif1}
\end{equation}
and the next players' round duration conditions
(\ref{stop2a}, \ref{stop2b}) become
\begin{equation}
    \forall(j\in\mathbb{N},j<J_{1},\vec{d}_{j}\in\partial B^{+}(D)), 
    u(\ldots\vec{d}_{j})=-\infty.
    \label{uf2}
\end{equation}

\subsection{Optimal strategies \label{optstr}}
 
The greatest utility, drawn from any event-terminated history, is
\begin{eqnarray}
    \forall j\in\mathbb{N},u(\ldots\vec{d}_{j},\vec{d}_{j+1/2}) 
    &=&\max_{\vec{d}_{j+1}}
    u(\ldots\vec{d}_{j},\vec{d}_{j+1/2},\vec{d}_{j+1}),
    \label{u->^u1}\\
    u(.)&=&\max_{\vec{d}_{1}}u(.,\vec{d}_{1}).
    \label{u->^u}
\end{eqnarray}
The latter equation, where $\vec{d}_{1}$ is a dummy variable, is a 
more formal expression of the former.  The nature of the dummy 
variable is shown by its index (state for integer, event for 
half-integer).  The set of states, corresponding to optimal decisions, 
is
\begin{equation}
    S_{u}(.)\equiv\argmax_{\vec{d}_{1}}u(.,\vec{d}_{1}).
    \label{?S_u}
\end{equation}

A player's mixed strategy consists in choosing randomly between many 
decisions, according to a causal probability law,
\begin{equation}
    \vec{d}\mapsto P(.,\vec{d})\equiv\mathcal{P}(\vec{d}_{1}=\vec{d}|.),
    \sum_{\vec{d}_{1}} P(.,\vec{d}_{1})=1.
    \label{?P}
\end{equation}
$\mathcal{P}(X)$ means the probability of the event $X$.  The optimal 
mixed strategies are such that the support of the probability law 
(\ref{?P}) is a subset of $S_{u}(.)$ (among them are pure optimal 
strategies).

From the von Neumann-Morgenstern theorem,
\begin{equation} 
    u(.)=\sum_{\vec{d}_{1/2}}p(.,\vec{d}_{1/2})u(.,\vec{d}_{1/2}).
    \label{^u->u}
\end{equation}
where $p$ is a causal probability law, expressing the providential 
strategy and rules.  Because of utility conditions, such as 
(\ref{mortouvif1}, \ref{uf2}), there are, in (\ref{^u->u}), products 
$p\times u$ of the undetermined form $0\times\infty$, which ought to 
be replaced by zero (or the summation ought to be properly 
restricted).

Combining (\ref{u->^u}, \ref{^u->u}), or conversely,
\begin{eqnarray}
    u(.)&=&\max_{\vec{d}_{1}}\sum_{\vec{d}_{3/2}}
    p(.,\vec{d}_{3/2})u(.,\vec{d}_{1},\vec{d}_{3/2}),
    \label{^u->^u}\\
    u(.)&=&\sum_{\vec{d}_{1/2}}p(.,\vec{d}_{1/2})\max_{\vec{d}_{1}}
    u(.,\vec{d}_{1/2},\vec{d}_{1}).
    \label{u->u}  
\end{eqnarray}
The composition of $\max-\moy$ operations names the algorithm, which 
is the classical zero-sum game $\max-\min$, where the rational 
opponent has been replaced by neutral providence.  (\ref{^u->^u}, 
\ref{u->u}) are consistent with (\ref{udd=ud}): after the judgment, 
they simply repeat the utility forever, so that the $\max-\moy$ 
operations can be chained ad infinitum, no matter the end of round.  
Thus, the judgment can be arbitrarily postponed, without affecting 
strategy.  If judgment times have an upper bound (e.  g. the number of 
fates is finite), then all judgments can be postponed until a 
(collective) last judgment a time $J\in\mathbb{N},J\ge\max 
J_{\varphi}$, e.  g. the maximum round duration.

Relaxing the rule (\ref{d0=0}), and taking $\vec{d}_{0}$ as a 
parameter, the problem of fate management, i.  e. finding optimal 
strategies, is \emph{self-similar} under time-shifts, except for the 
parameter ``renormalization'' (as in statistical physics)
\begin{equation}
    (D_{0},J,\vec{d}_{0},(\ldots\vec{d}_{j}),(p,u)(.)) 
    \rightarrow(D_{j},J-j,\vec{d}_{j},(),(p,u)(\ldots\vec{d}_{j},.)).
    \label{renorm}
\end{equation}

Let $\chi_{V}$ be the characteristic function of $V\subset B(D)$.  The 
round providential strategy is determined by (\ref{multinomial}) and
\begin{equation}
    p(.,\vec{d}_{0},\vec{d}_{1/2})\equiv p(\vec{d}_{1/2})
    \chi_{\partial B^{+}(D_{0})}(\vec{d}_{1/2}).
    \label{?p}
\end{equation}
The expected utility is computed with (\ref{u->^u}, \ref{^u->u}), from 
the last judgment backward in time:
\begin{equation}
    u(\ldots\vec{d}_{J-1},\vec{d}_{J-1/2}),
    u(\ldots\vec{d}_{J-1})\ldots
    u(\vec{d}_{0},\vec{d}_{1/2},\vec{d}_{1}), 
    u(\vec{d}_{0},\vec{d}_{1/2}),
    u(\vec{d}_{0}),
\end{equation}
e. g., for $J=3$, and using (\ref{?p}),
\begin{eqnarray}
    u(\ldots\vec{d}_{3/2})
    =\max_{\vec{d}_{2}}&&\sum_{\vec{d}_{5/2}\in\partial B^{+}(D_{2})}
    p(\vec{d}_{5/2})u(\ldots\vec{d}_{3/2},\vec{d}_{2},\vec{d}_{5/2}),
     \nonumber
\\
    u(\vec{d}_{0},\vec{d}_{1/2})
    =\max_{\vec{d}_{1}}&&\sum_{\vec{d}_{3/2}\in\partial B^{+}(D_{1})}
    p(\vec{d}_{3/2})
    u(\vec{d}_{0},\vec{d}_{1/2},\vec{d}_{1},\vec{d}_{3/2}).
    \nonumber
\\
    u(\vec{d}_{0})=&&\sum_{\vec{d}_{1/2}\in\partial B^{+}(D_{0})}
    p(\vec{d}_{1/2})u(\vec{d}_{0},\vec{d}_{1/2}).
    \nonumber
\end{eqnarray}

\section{Fate as a stochastic process\label{stopro}}

For a given strategy, what is the presence density (of the die 
system in a subset of phase space)?  What is the expectation of an 
arbitrary utility, for which the given strategy is not necessarily 
optimal?  

\subsection{The Kolmogorov equation on expected utility}

Fate is a stochastic process, not only because it contains random 
events (the probability law $p$), but also random decisions, according 
to mixed strategies (the probability law $P$).  For any causal process 
like (\ref{destin}), the sequence of histories
\begin{equation}
    (\vec{d}_{0}),
    (\vec{d}_{0},\vec{d}_{1/2}),
    (\vec{d}_{0},\vec{d}_{1/2},\vec{d}_{1})\ldots
    \nonumber
\end{equation}
is a discrete Markov chain, for which classical results are available 
\cite[ch.  6]{Parzen}, \cite[ch.  15]{Papoulis}, originating mostly 
from Brownian motion studies \cite[ch.  15]{Reif}.

The fate stochastic evolution equation, the Langevin equation, is just 
a random sum, obeying (\ref{toutrejouer}, \ref{mortouvif}):
\begin{eqnarray*}
    \vec{d}_{j+1}&=&\vec{d}_{j}+\hat{\vec{d}}_{j+1},
    \\
    \mathcal{P}(\hat{\vec{d}}_{j+1}=\vec{d}|\ldots\vec{d}_{j},\vec{d}_{j+1/2})
    &=&P(\ldots\vec{d}_{j},\vec{d}_{j+1/2},\vec{d}_{j}+\vec{d}),
    \\
    \mathcal{P}(\vec{d}_{j+1/2}=\vec{d}|\ldots\vec{d}_{j})
    &=&p(\ldots\vec{d}_{j},\vec{d}).
\end{eqnarray*}
$\hat{\vec{d}}_{j+1}$ is a random source term, conditioned by history, 
according to the mixed strategies $P,p$.  $\vec{d}_{j}$ undergoes a 
strategy-driven Brownian motion as, for example, a charged Brownian 
particle driven by electrophoresis.

The Chapman-Kolmogorov equation yields the probability of transition, or 
jump, from one state to the other, in one time step:
\begin{equation}
    \sigma(.,\vec{d}_{0}\curvearrowright\vec{d}_{1}) 
    \equiv\mathcal{P}(\vec{d}_{1}=\vec{d}|.,\vec{d}_{0}) 
    =\sum_{\vec{d}_{1/2}} 
    p(.,\vec{d}_{0},\vec{d}_{1/2})P(.,\vec{d}_{0},\vec{d}_{1/2},\vec{d}).
    \label{CK}
\end{equation}

Let $P$ be a player's mixed strategy, possibly not optimal.  From the 
von Neumann-Morgenstern theorem, twice applied,
\begin{equation}
    u(.,\vec{d}_{0})=\sum_{\vec{d}_{1/2}}p(.,\vec{d}_{0},\vec{d}_{1/2})
    \sum_{\vec{d}_{1}}P(.,\vec{d}_{0},\vec{d}_{1/2},\vec{d}_{1})
    u(.,\vec{d}_{0},\vec{d}_{1/2},\vec{d}_{1}).
    \label{tu0}
\end{equation}

Reversing the order of summation, using (\ref{CK}) and assuming that 
\emph{utility does not depend on events, but only on states}, which is 
true in the 421 round, I obtain the Kolmogorov equation on the 
expected utility:
\begin{equation}
    u(.,\vec{d}_{0})
    =\sum_{\vec{d}_{1}}\sigma(.,\vec{d}_{0}\curvearrowright\vec{d}_{1})
    u(.,\vec{d}_{0},*,\vec{d}_{1}).
    \label{tu}
\end{equation}
(By hypothesis, $u$ does not depend on $*$.)

As opposed to the $\max-\moy$ algorithm, (\ref{tu}) does not produce 
any decision, but, given the mixed strategies $P,p$ (effective through 
$\sigma$), determines the expectation of any utility, for which $P$ 
may not be optimal.

Nevertheless, if $P$ is optimal, from (\ref{u->^u}) and (\ref{?P}), 
there is an equality, between operators on 
$u(.,\vec{d}_{0},\vec{d}_{1/2},\vec{d}_{1})$:
\begin{equation}
    \max_{\vec{d}_{1}}
    =\sum_{\vec{d}_{1}}P(.,\vec{d}_{0},\vec{d}_{1/2},\vec{d}_{1}).
    \label{max=sum}
\end{equation}
Taking (\ref{max=sum}) into (\ref{tu0}) returns (\ref{u->u}).

\subsection{The Fokker-Planck equation on presence density}

I define the state fate $\psi\equiv(\vec{d}_{j})_{j=0,1\ldots}$ (fate 
with only states, not events).  From (\ref{CK}),
\begin{equation}
    \mathcal{P}(\psi=(.,\vec{d}_{0},\vec{d}_{1}))
    =\sigma(.,\vec{d}_{0}\curvearrowright\vec{d}_{1})
    \mathcal{P}(\psi=(.,\vec{d}_{0})),
    \label{CK1}
\end{equation}
so that the sequence of past states
\begin{equation}
    (\vec{d}_{0}),
    (\vec{d}_{0},\vec{d}_{1}),
    (\vec{d}_{0},\vec{d}_{1},\vec{d}_{2})\ldots
    \nonumber
\end{equation}
also is a Markov chain.

Summing (\ref{CK1}) over all state fates converging to the same state 
$\vec{d}$ at time $j+1$ gives the presence density 
$\rho_{j+1}(\vec{d})$:
\begin{eqnarray}
    \rho_{0}(\vec{d})&=&\delta_{\vec{d},\vec{d}_{0}},
    \label{rho0=}\\
    \forall j\in\mathbb{N},
    \rho_{j+1}(\vec{d})&=&\sum_{\vec{d}_{0}\ldots\vec{d}_{j}} 
    \mathcal{P}(\psi=(\vec{d}_{0}\ldots\vec{d}_{j})).
    \label{rhoj=}
\end{eqnarray}
In the round, from (\ref{asymp}), $\rho_{j}$ is stationary, as soon as 
$j\ge J$.

I assume that \emph{utility is a function of state and time only}, 
less general than causal (\ref{causale}):
\begin{equation}
    u(\ldots\vec{d}_{J_{\varphi}})=u_{J_{\varphi}}(\vec{d}_{J_{\varphi}}).
    \label{u(t,d)}
\end{equation}
The end-of-round utility is indeed of the kind (\ref{u(t,d)}), because 
end-of-set ranking (see the rules) only depends on round results, not 
on intermediary states and, for the first player, the effective round 
duration.

For all player's optimal (or rational) mixed strategy $P$ derived 
from a utility of the kind (\ref{u(t,d)}),
\begin{eqnarray}
    P_{j}(\vec{d}_{j},\vec{d}_{j+1/2},\vec{d}_{j+1})
    &\equiv&P(.,\vec{d}_{j},\vec{d}_{j+1/2},\vec{d}_{j+1}),
    \label{?Pj}\\
    \sigma_{j}(\vec{d}_{j}\curvearrowright\vec{d}_{j+1})
    &\equiv&\sigma(.,\vec{d}_{j}\curvearrowright\vec{d}_{j+1}) 
    \label{?sigmaj}
\end{eqnarray} 
Taking (\ref{u(t,d)}, \ref{?sigmaj}) into (\ref{tu}) allows to extend 
(\ref{u(t,d)}) to all time (for the expected utility), by induction:
\begin{equation}
    \forall j\in\mathbb{N},u_{j}(\vec{d}_{j})\equiv u(\ldots\vec{d}_{j}).
    \label{uj(d)}
\end{equation}
The process $(j,\vec{d}_{j},\vec{d}_{j+1/2})$ is Markovian.

The consequence (\ref{?sigmaj}) of (\ref{u(t,d)}), taken into 
(\ref{rhoj=}), allows to express $\rho_{j+1}(\vec{d}_{j+1})$ as a 
functional on $\rho_{j}$:
\begin{equation}
    \rho_{j+1}(\vec{d}_{j+1})=\sum_{\vec{d}_{j}}
    \rho_{j}(\vec{d}_{j})\sigma_{j}(\vec{d}_{j}\curvearrowright\vec{d}_{j+1}),
    \label{trho} 
\end{equation}
the Fokker-Planck equation.

As opposed to (\ref{trho}), (\ref{tu}) does not need (\ref{u(t,d)}).  
Nevertheless, with (\ref{u(t,d)}), (\ref{tu}), becomes
\begin{equation}
    u_{j}(\vec{d}_{j}) 
    =\sum_{\vec{d}_{j+1}}\sigma_{j}(\vec{d}_{j}\curvearrowright\vec{d}_{j+1}) 
    u_{j+1}(\vec{d}_{j+1}),
    \label{tu1}
\end{equation}
adjoint to (\ref{trho}).

(\ref{trho}, \ref{tu1}) are the evolution equations, adjoint to each 
other, linear, unstationary, of presence density and expected utility.  
Their inputs are a player's mixed strategy and utility.

\subsection{Computing result probabilities by duality\label{dual}}

Let $\mathcal{F}\equiv\mathcal{F}(B^{+}(D),\mathbb{Q})$ be the space 
of numerical functions on $B^{+}(D)$, with the scalar product
\begin{equation}
    \forall(f,g)\in\mathcal{F},
    \<f,g\>\equiv\sum_{\vec{d}\in B^{+}(D)}f(\vec{d})g(\vec{d}).
    \label{?<,>}
\end{equation}
$\sigma_{j}$ is an operator, a linear endomorphism on $\mathcal{F}$, 
fully determined by the Markovian matrix 
$\sigma_{j}(\vec{d}_{0}\curvearrowright\vec{d}_{1})$.  Its transposed operator is 
$\sigma_{j}^{t}$, of matrix \[\sigma_{j}^{t}(\vec{d}_{1}\curvearrowright\vec{d}_{0}) 
\equiv\sigma_{j}(\vec{d}_{0}\curvearrowright\vec{d}_{1}).\] In operator notation, 
(\ref{trho}, \ref{tu1}) become
\begin{equation}
    \rho_{j+1}=\sigma_{j}^{t}\rho_{j},u_{j}=\sigma_{j}u_{j+1}.
    \nonumber
\end{equation}

As $\sigma_{j}^{t}$ et $\sigma_{j}$ are adjoint to each other, the 
expected utility follows a conservation law:
\begin{eqnarray}
    \<u_{j},\rho_{j}\>=\<\sigma_{j}u_{j+1},\rho_{j}\>&=&
    \<u_{j+1},\sigma_{j}^{t}\rho_{j}\>=\<u_{j+1},\rho_{j+1}\>,
    \nonumber\\
    \<u_{j},\rho_{j}\>
    &=&\<u_{0},\rho_{0}\>=u_{0}(\vec{d}_{0}).
    \label{urho=u0}
\end{eqnarray}
The last equality is a consequence of (\ref{rho0=}).  Given the 
player's mixed strategy $P$, (\ref{urho=u0}) holds for any utility.

The direct computation of $\<u_{j},\rho_{j}\>$ consists in solving for 
$\rho_{j}(\vec{d}_{j})$ the Fokker-Planck equation, which must be 
repeated, to complete the scalar product, at least for all 
$\vec{d}_{j}$ where $u_{j}$ does not vanish.  More shrewdly, 
$\<u_{j},\rho_{j}\>$ can be computed indirectly, as the r.  h. s.  of 
(\ref{urho=u0}): the Kolmogorov equation is solved only once for the 
expected utility at the trunk of the fate tree, or the initial 
expected utility.  The indirect computation is faster than the direct 
computation, by a factor which is the cardinal of the support of 
$u_{j}$.  The indirect computation benefits from the unicity of the 
fate tree, and the diffusive growth of the support of $u_{j}$.

Moreover, to obtain the Kolmogorov algorithm from the $\max-\moy$ 
algorithm, one merely has to replace, in (\ref{u->u}), the operator 
$\max$ appearing at the l.  h. s.  of (\ref{max=sum}), by the operator 
$\str$ appearing at the r.  h. s.  of (\ref{max=sum}).  (These 
operators differ if $P$ is not optimal.)  The Kolmogorov equation is 
thus solved by a $\str-\moy$ algorithm.

Here are examples of using the Kolmogorov equation and 
(\ref{urho=u0}):
\begin{enumerate}
    \item The probability of the result to be in $V\subset B^{+}(D)$ 
    (independently of time) is the initial expectation of the 
    stationary utility $u_{j}=\chi_{V}$.

    \item The probability of $D_{j_{0}}$ is the initial expectation of 
    the utility $u_{j}=\delta_{j,j_{0}}\chi_{\partial 
    B^{+}(D-D_{j_{0}})}$.
\end{enumerate}

\subsection{Analogy with linear transport theory}

The round is a linear transport phenomenon, with respect to the face 
variable.  Face, expected utility, presence density, transition 
probability correspond respectively, in transport theory \cite{Case}, 
to phase (position, velocity), importance \cite{Lewins1}, flux and 
cross section.  Harris \cite{Harris} shows that a monokinetic particle 
population grown by branching (e.  g. neutrons produced by nuclear 
fission) follows a Galton-Watson process.  Similarly, in appendix 
\ref{GW}, I discuss the Galton-Watson character of the first player's 
live dice population $D_{j}$.

\section{Simple optimal policies for one-goal utilities\label{1goal}}

Taking for goal a unique combination $\vec{d}^{*}\in\partial 
B^{+}(D)$, the utility is a binary Kronecker function 
$\delta_{\vec{d}^{*},.}$ (modulo an affine transform), and optimal 
strategies are simply constructed.

\subsection{The ratchet and Bernoulli policies\label{ss:cliquet}}

I examine two first player's policies, with a one-goal utility:
\begin{enumerate}
    \item The Bernoulli policy consists in accumulating no die, unless 
    the goal has been attained (then, all dice are accumulated); the 
    cast sequence is a stationary Bernoulli process (a sequence of 
    independent trials terminated by success or failure).

    \item The ratchet policy consists in putting aside as many dice as 
    possible, contributing to the goal:
\begin{eqnarray}
    \forall(j\in\mathbb{N},j+1<J), 
    \vec{d}_{j+1}&=&\vec{d}^{*}\wedge(\vec{d}_{j}+\vec{d}_{j+1/2}), 
    \label{cliquet}\\
    P_{j}(\vec{d}_{j},\vec{d}_{j+1/2},\vec{d}_{j+1})
    &=&\delta_{\vec{d}_{j+1},\vec{d}^{*}\wedge(\vec{d}_{j}+\vec{d}_{j+1/2})}.
    \label{Pcliquet}
\end{eqnarray}
\end{enumerate}
($\wedge$ is the infix notation of the minimum in the partially 
ordered space $\mathbb{Z}^{F}$, generalizing, in Lagrangian notation, 
the ensemble intersection $\cap$.)

The ratchet strategy towards $\vec{d}^{*}$ is optimal, with respect to 
the $\vec{d}^{*}$-goal utility, if and only if $p$ decreases in 
$B^{+}(D)$.  This means that as many dice as possible should be 
accumulated, in order to maximize the success probability at any 
future time.  For unloaded dice, from (\ref{multinomial}),
\begin{equation}
    \forall(\vec{d}\in B^{+}(F),\vec{d}+\vec{e}_{1}\in B^{+}(F)),
    \frac{p(\vec{d}+\vec{e}_{1})}{p(\vec{d})}
    =\frac{1}{F}\frac{|\vec{d}|+1}{d_{1}+1}\le 1,
    \label{pdecroi}
\end{equation}
i. e.  $p$ decreases on $B^{+}(F)$.  The ratchet strategy is optimal 
if and only if $D\le F$, strictly if and only if $D<F$.

For example, with $D=3<F=6,J>1,\vec{d}^{*}=421,\vec{d}_{1/2}=651$, the 
ratchet decision (to accumulate 1) is optimal, because $p(421)<p(42)$ 
(it will be easier to obtain 42 than 421).  With 
$D=3>F=2,\vec{d}^{*}=211,\vec{d}_{1/2}=222$, the Bernoulli decision (to 
replay all dice) is optimal, because 
$p(11)=1/F^{2}=2/8<p(211)=3/F^{3}=3/8$.  With 
$D=F=2,\vec{d}^{*}=21,\vec{d}_{1/2}=11$, both Bernoulli and ratchet 
decisions are optimal.

A next player's maximum round duration is imposed.  In case of a 
premature success, he is in a \emph{dilemma}, having to decide between 
equally unpleasant ways of breaking the goal, obtained too early.  For 
$F>2$, optimal decisions consist in replaying any one die; the number 
of pure optimal strategies is thus the number of distinct faces in the 
goal combination, at the power $J-1$.  If the goal is a brelan, then 
no dilemma exists.

\subsection{Optimal one-goal strategy result probabilities 
\label{1goalprob}}

For any strategy, I consider the probability to obtain any result, e.  
g. 111 after three casts.  According to section \ref{dual}, this 
probability is the initial expected utility, determined by the 
Kolmogorov equation and the final condition of a Kronecker utility on 
the result.  This probability depends on the player $i=1,2$ (first or 
next), the (renormalized) maximum round duration $J_{1}$, the player's 
mixed strategy $P$, the delay $j$, and the result $\vec{d}$:
\begin{equation}
    p_{i}(J_{1},P,j,\vec{d}),0\le j\le J_{1}\le J,\vec{d}\in B^{+}(D).
    \label{?p()}
\end{equation}

The set of result probabilities, for all possible pure strategies and 
$(D,F,J)=(3,6,3)$, is (much larger than the fate tree, itself very 
large and) too large to be extensively listed.  Thus, I will work on a 
reduced strategy subset, for which a reasonable choice is the set of 
optimal one-goal strategies, for all possible goals.  As far as the 
goal determines the optimal strategy, the variable $P$ in (\ref{?p()}) 
is simply replaced by the goal $\vec{d}^{*}$:
\begin{equation}
    p_{i}(J_{1},\vec{d}^{*},j,\vec{d}), 0\le j\le J_{1}\le J,
    (\vec{d},\vec{d}^{*})\in B^{+}(D)^{2}
    \label{?p()1}
\end{equation}  
which looks like the Markovian matrix of section \ref{dual}, except 
that $\vec{d}^{*}$ is not actual, but contemplated.  There are 
diagonal ($\vec{d}=\vec{d}^{*}$) and non-diagonal result 
probabilities.

For the first player, the optimal one-goal strategy is unequivocally 
defined by the goal ($D<F$: the ratchet) and the function $p_{1}$ is 
defined everywhere.  This in not true for $p_{2}$, because of 
dilemmas.  However, next player diagonal probabilities are unaffected 
by dilemmas, so that $p_{2}$ is defined on the diagonal, 
$\vec{d}=\vec{d}^{*}$; it is even defined for all 
$(\vec{d}^{*},\vec{d})$, if and only if $\vec{d^{*}}$ is a brelan, 
since brelans do not produce dilemma, as noticed at end of section 
\ref{ss:cliquet}.

Here are a few properties of the functions $p_{i}$:
\begin{eqnarray}
    p_{i}(0,\vec{d}^{*},0,\vec{d})&=&\delta_{\vec{d}^{*},\vec{d}},
    \label{propp}\\
    p_{i}(J,\vec{0},j,\vec{0})&=&\delta_{j,0},
    \label{p_i(J0j0)}\\
    p_{i}(1,\vec{d}^{*},1,\vec{d})&=&p(\vec{d}),
    \label{p_i(11)}\\
    p_{i}(J,\vec{d}^{*},j,\vec{d})&=&0,j<J,\vec{d}^{*}\neq\vec{d},
    \nonumber\\
    p_{1}(J,\vec{d},j,\vec{d})&=&p_{i}(j,\vec{d},j,\vec{d}),j<J,
    \nonumber\\
    p_{2}(J,\vec{d}^{*},j,\vec{d})&=&0,j<J,
    \label{p_2=0}\\
    \sum_{\vec{d}\in\partial B^{+}(|\vec{d}^{*}|)}
    \sum_{j=0}^{J}p_{1}(J,\vec{d}^{*},j,\vec{d})&=&1.
    \nonumber\\
    \sum_{\vec{d}\in\partial B^{+}(|\vec{d}^{*}|)}
    p_{2}(j,D\vec{e}_{f},j,\vec{d})&=&1.    
    \nonumber
\end{eqnarray}

Let the cumulative diagonal probability be
\begin{equation}
    s_{i}(J,\vec{d})\equiv\sum_{j=1}^{J}p_{i}(J,\vec{d},j,\vec{d}).
    \label{?pcumul}
\end{equation}
Because of the next players' round duration condition 
\begin{equation}
    \forall J>1,s_{1}(J,\vec{d})>s_{2}(J,\vec{d})=p_{2}(J,\vec{d},J,\vec{d})
    >p_{1}(J,\vec{d},J,\vec{d}).  
    \nonumber
\end{equation}

To reduce the $p_{i}$ computational domain, I use invariance with 
respect to face permutations (for unloaded dice).  Firstly, diagonal 
probabilities depend on only one combination.  As in 
(\ref{multinomial}), two combinations are equivalent, modulo the 
functions $\vec{d}\mapsto p_{i}(J,\vec{d},j,\vec{d})$, for all 
$(i,J,j)$, if and only if their occupation numbers (Lagrangian 
components) form the same combination, e.  g. $441\sim 655$.  With 
$(D,F)=(3,6)$, the quotient set contains three classes: that of 
brelans ($\ni 111$), that of sequences ($\ni 123$)\footnote{I do not 
mean that all combination in the class of sequences is a sequence.}, 
that of pairs ($\ni 112$).  Secondly, non-diagonal probabilities 
depend on a couple of combinations.  Two \emph{couples} of 
combinations are equivalent, modulo the functions 
$(\vec{d}^{*},\vec{d})\mapsto p_{i}(J,\vec{d}^{*},j,\vec{d})$, for all 
$(i,J,j)$, if and only if their couples of occupation numbers form the 
same combination, e.  g. $(421,442)\sim (321,211)$.  A face 
permutation transforms a next player's optimal one-goal strategy into 
another, possibly different if the goal is not a brelan.

Taking into account (\ref{propp}) and face permutation invariance, the 
result probabilities (\ref{?p()1}) are computed, for 
$(D,F,J)=(3,6,3)$, by applying $\str-\moy$ on optimal 
$\vec{d}^{*}$-goal strategies and $\vec{d}$-Kronecker utilities.  As a 
consequence of self-similarity (\ref{renorm}), the probabilities after 
the initial time ($J_{1}<J$), are obtained as intermediary results in 
the computation of a priori probabilities ($J_{1}=J$).  The results 
are presented in the probability charts \ref{diaPro}, \ref{nonDiaPro}, 
\ref{nonDiaPro1}, \ref{nonDiaPro2}, \ref{nonDiaPro3} (appendix 
\ref{mma}), which do not fill more than a few pages thanks to the 
extensive use of face permutation invariance and other properties 
(\ref{propp}\ldots).  There are 31 classes of three-die combination 
couples (including the three diagonal classes).

\section{Goal identification programming}

I will propose heuristic policies, based on the global maximization of 
expected utility, with respect to the subset of optimal one-goal 
strategies, for which result probabilities were obtained in the last 
section.

\subsection{Motivation: bounded complexity}

The $\max-\moy$ backward induction algorithm is optimal, short, but 
the number of numerical operations per time step, already large for 
$(D,F,J)=(3,6,3)$, is unbounded as a function of the maximum round 
duration $J$.  Information theory \cite{Brillouin,Delahaye} teaches 
that a message will be transmitted faster by a specialized code.  
$\max-\moy$ backward induction is slow, for the general reasons that 
it is unspecialized (and optimal).

To speed-up policy, possibly at the expense of brevity and optimality, 
specialization is necessary.  For example, consider the game of Nim 
\cite[§ 1.3]{Dutta}: besides $\max-\moy$ backward induction, a 
stratagem is found, based on congruence, producing optimal strategies, 
with a bounded number of operations per time step.  The ratchet 
($D<F$) would be a stratagem of 421, if only the goal were known.

I propose to identify the goal, rigorously, by considering not only 
the utility, but also the result probabilities (\ref{?p()1}), obtained 
in section \ref{1goal}.  I will obtain goal identification heuristic 
policies, that may be considered as quasi-Markovian, from the remark 
following (\ref{?p()1}).  Roughly, they transfer the complexity of 
$\max-\moy$ backward induction to the result probabilities, with the 
advantage that the latter can be compiled once for all (and the 
inconvenience that they must be remembered).

For a one-goal utility, goal identification is simple.  For a constant 
utility, as well: any goal is optimal.  Difficulties are thus with 
utilities somewhere between peaked and flat, ``\emph{fuzzy}'', e.  g. 
with peaks of about the same height, playing the roles of attractors, 
that one has to choose between.\footnote{Like Buridan's donkey, 
starving from hesitating between bushels of oats and water.}

\subsection{Reduced horizon\label{choix}}

I consider a time and state dependent utility, as in (\ref{u(t,d)}, 
\ref{uj(d)}), in a round of maximum duration $J$.  $\vec{d}_{0}$ 
is the state at time $j_{0}\le J$.  I define the ``evaluation 
function'',
\begin{equation}
    u^{*0}_{j_{0}}(\vec{d}_{0})\equiv\max_{\vec{d}^{*}\in\partial B^{+}(D_{0})} 
    \sum_{j=0}^{J-j_{0}}p_{i}(J-j_{0},\vec{d}^{*},j,\vec{d}^{*})
    u_{j_{0}+j}(\vec{d_{0}}+\vec{d}^{*}),
    \label{?u^*}
\end{equation}
where $j$ is the renormalized time and $D_{0}=D-|\vec{d_{0}}|$.  Evaluation 
functions are often used in stage game (chess, othello, 
checkers\ldots) programming, but they are usually defined empirically, 
unlike (\ref{?u^*}), which is probabilistic.

To take into account serendipity -- that a result other than the goal 
may be not so bad, after all -- (\ref{?u^*}) is improved:
\begin{equation}
    u^{*1}_{j_{0}}(\vec{d}_{0})\equiv\max_{\vec{d}^{*}\in\partial 
    B^{+}(D_{0})} \sum_{j=0}^{J-j_{0}}\sum_{\vec{d}\in\partial 
    B^{+}(D_{0})} p_{1}(J-j_{0},\vec{d}^{*},j,\vec{d}) 
    u_{j_{0}+j}(\vec{d_{0}}+\vec{d}),
    \label{?u^*s}
\end{equation}
which cannot be used for next players, because of dilemmas.  For all 
$\vec{d}_{0}\in\partial B^{+}(D)$, considering (\ref{p_i(J0j0)}), the 
evaluation functions (\ref{?u^*}, \ref{?u^*s}) simply return the 
utility.

$\max-\moy$ backward induction is particularly slow, because it needs 
to completely analyze the round even before its first decision.  Hence 
the idea that short-sighted policies may be faster.  At time 
$j_{0}\in\mathbb{N}$, a horizon $h\in\mathbb{N}$ may be chosen, such 
that $j_{1}=j_{0}+h\le J$, and the round is virtually terminated at 
$j_{1}$, taking for ersatz utility the evaluation function 
$u^{*s}_{j_{1}}$ given by (\ref{?u^*}) or (\ref{?u^*s}), depending on 
the serendipity bit $s\in\{0,1\}$.  With $j_{1}=J-1$, considering 
(\ref{p_i(11)}), (\ref{?u^*s}) reproduces the deepest $\max-\moy$ 
iteration, so that an optimal strategy is generated.

I will further examine $h=0,1$.  With $h=0$, the goal is found by 
maximizing $u^{*s}_{j_{0}}$, independently of the first event.  With 
$h=1$, as there is no interest in thinking before casting the dice, 
the decision $\vec{d}_{1}$ is rather taken after the first event 
$\vec{d}_{1/2}$, according to
\begin{equation}
    \max_{\vec{d}_{1}}
    u^{*s}_{j_{0}+1}(\vec{d}_{1}).
    \label{maxd_1}
\end{equation}
In case of many optimal decisions in (\ref{maxd_1}), the corresponding 
states, written as increasing Lagrangian lists, e.  g. 124, are 
discriminated according to the lexicographic order (only pure 
strategies are generated).  In case of many optimal goals in 
(\ref{?u^*}) or (\ref{?u^*s}), we need not discriminate between 
them, and the policy reproduces the human character of 
\emph{duplicity}.  Dilemma implies duplicity, but the converse is 
false.

\subsection{Dynamic programming and goal revision}

The strategy may be revised to take into account new events, which is 
an instance of dynamic programming \cite{Bellman} or belief revision 
\cite{Fabiani}, realizing a feedback of fate on strategy.  By 
self-similarity of the round, a policy may be applied at any time, 
with suitable parameter renormalization.  Self-similar revision based 
on the $\max-\moy$ backward induction policy would just confirm the 
optimal strategy, computed a priori: it is therefore useless.  Only 
fallible policies are worth revising.

A heuristic policy of horizon $h\ge 1$ forecasts, at any given time, 
only the next $h$ decisions.  Thus, it must be run with the period at 
least $h$.  The revised serendipitous goal identification policy of 
horizon $h$ is optimal in its last $h$ decisions.  The goal 
identification policy with $h=0$ does not require revision and is very 
simple (short and fast).  It may be the only rational policy, simple 
enough for unaided human players in normal game conditions.

\subsection{Policy benchmark and interpretation\label{opt}}

For $(D,F,J)=(3,6,3)$, I consider a few increasingly fuzzy stationary utilities:
\begin{enumerate}
    \item $u=\delta_{123}$, a one-goal utility,
    
    \item $u=\delta_{123}+\delta_{224}+\delta_{345}$, a three-goal 
    utility,
    
    \item $u=t$, the transfer function defined by table \ref{transf} in 
    appendix,

    \item the sum of faces.
\end{enumerate}
These utilities are unrealistic, in the sense that they may not be 
possible within a real 421 set (see section \ref{uti}).  I consider 
the policies: $\max-\moy$ backward induction, and the four goal 
identification policies $(h,s)\in\{0,1\}^{2}$; the $h=0$ policies are 
without revision.

From the final utility, on the leaves of the fate tree, every policy 
yields a pure strategy, and its initial expected utility $u_{0}$, on 
the trunk, is obtained by solving the Kolmogorov equation exactly, 
with the $\str-\moy$ algorithm.  Optimality is defined as the ratio of 
the expected utility, over the first player optimal expected utility 
$u_{0r}$.  The numerical results (approximated by decimal numbers) are 
copied from \cite{421mma} into the tables \ref{u=1obj}, \ref{u=3obj}, 
\ref{u=t}, \ref{u=somf}.

\newcommand{\opttab}[6]{
\vspace{1em}
\begin{tabular}{|c|c|c|c|c|}   
    \hline\multicolumn{2}{|c|}{$u_{0r}=#1$}&\multicolumn{2}{|c|}{$u_{0}/u_{0r}$}\\
    \hline\multicolumn{2}{|c|}{policy}&\multicolumn{2}{|c|}{player}\\
    \hline 
    horizon&serendip. &first& next\\\hline
    0         & 0            & #3            \\\hline
    0         & 1            & #4 &         \\\hline
    1         & 0            & #5            \\\hline
    1         & 1            & #6 &         \\\hline
    \multicolumn{2}{|c|}{$\max-\moy$}
                                & 1   & #2    \\\hline
\end{tabular}
}
\begin{table}[tbp]
    \centering
    \caption{123 one-goal utility}
    \label{u=1obj} 
    \opttab{0.22811}{0.57858}{1&0.57858}{1}{1&0.57858}{1}
\end{table}
\begin{table}[tbp]   
    \centering
    \caption{123, 224, 345 three-goal utility}
    \label{u=3obj} 
    \opttab{0.32805}{0.49152}{0.73037&0.43734}{0.73037}
    {0.97777&0.47746}{0.98657}
\end{table}
\begin{table}[tbp]
    \centering
    \caption{utility = transfer function}
    \label{u=t}
    \opttab{3.7467}{0.77663}{0.90834&0.68812}{0.90834}
    {0.87962&0.68991}{0.99634}
\end{table}
\begin{table}[tbp]
    \centering
    \caption{utility = sum of faces}
    \label{u=somf}
    \opttab{14}{0.97321}{0.94194&0.92599}{0.96418}{0.75&0.85875}{0.99900}
\end{table}

Table \ref{u=1obj} confirms that for a one-goal utility, all goal 
identification policies are by definition optimal.  Compared to the 
first player, next players are handicapped, but less with a fuzzier 
utility.
The numerical results show a positive contribution of serendipity, 
much greater with the greater horizon and revision.  The contribution 
of horizon and revision is positive with serendipity.  Without 
serendipity, the contribution of horizon and revision is positive for 
peaked utilities, negative for fuzzy utilities (\ref{u=t}, 
\ref{u=somf}).

I take advantage of this effect to give a (less fuzzy) definition of 
fuzziness: a utility is fuzzy if and only if introducing horizon and 
revision without serendipity contributes negatively to its 
expectation.  Thus, I have constructed fuzzy utilities, for which 
introducing horizon and revision decreases the expected utility, even 
though it is more complex.  The response of expected utility with 
respect to complexity is non-increasing (this effect compares, in 
electricity, with a negative resistance).

\section{Conclusions}

The mathematics of fate in 421 leave as the only unsolved difficulty 
``bifurcations'', that maximizing the expected utility does not always 
determine a unique decision, as in next players' dilemmas.  Here is a 
toy example: a game with three players, P, A, B. If P says white, then 
A gives one euro to B; if P says black, then B gives one euro to A. P 
earns nothing anyway; A, B take no decision.  Maximizing P's expected 
utility does not determine its decision.  Introducing a mixed strategy 
amounts to consider P as a random generator, with unknown 
probabilities.  A classic postulate of statistical theories is to 
maximize the entropy or missing information \cite{Reif,Brillouin}, 
which here sets the probabilities of either outcome to $1/2$.  Are the 
postulates of mixed strategy and maximum entropy so easily acceptable?  
We cannot exclude hidden determinism or bias in P. For example, P may 
always choose the first answer in the lexicographic order (black), or 
P may have a secret agreement with A to share his gain.

Bounded complexity, similar to \emph{bounded rationality} in 
\cite{Walliser}, motivates heuristic policies, where characters close 
to actual human behavior are found, in agreement with \cite{Goeree}.  
These characters are fate, dilemma, goal identification and revision, 
restricted horizon, serendipity, duplicity and panic.  When the policy 
belongs to an organization, we are in management.  When an individual 
decides for himself, we are in psychology.  For example, the same 
mathematical effect is behind counterproductive management or panic.

Goal identification consumes a bounded number of operations per time 
step, whatever the round duration, because it does not resolve all 
decisions in the fate tree, but only those which are compatible with 
the present state, and before the horizon.  Goal identification is not 
generally optimal, as opposed to a common assertion in business 
courses.  Only $\max-\moy$ backward induction, which has no goal, just 
like random playing, is generally optimal.  In the round, the ratchet 
stratagem allows the immediate translation from goal to decision.  I 
used probability theory as the logic of goal identification, \`a la 
Jaynes \cite{Jaynes}.  Complexity hides in the result probabilities, 
to be compiled before playing, as a kind of training.

Depending on complexity resources and utility, policies may be 
variably applicable or good.  Starting from a given policy, one may 
increase optimality, by modifying its characters or the utility: this 
is the task of human resources management, when the policy is that of 
an individual taking decisions for a company, a manager.  The 
short-sighted manager ($h=0$) gets hardly any help from serendipity.  
The unserendipitous manager should avoid fuzzy utilities and favor 
precise goal assignments.  I obtain examples of counterproductive 
management: with a fuzzy utility and no serendipity, goal revision 
dramatically reduces the optimality.  The role of serendipity was 
pointed out, on purely qualitative ground, by N. Wiener, about 
scientific and technical invention \cite{Wiener}.  The present work 
also pertains to Wiener's cybernetics.

Rationality can be further reduced.  At the extreme, the fool manager 
can be trusted only for a flat utility.  The study of irrational or 
illogical but actual behavior is the task of sophistry 
\cite{Schopenhauer}.  It may be quite useful in game practice, to 
produce best responses.

I thank researchers of the GREQAM in Marseilles, for fruitful 
discussions.

\appendix

\section{The (tentative) rules of 421\label{reg}}

I define the game, from oral tradition and \cite{valet,jeujura}.  The 
hardware consists of three dice and eleven tokens, initially in a pot.  
There are two or more players who can always see the positions of dice 
and tokens.

In the first part of the game, the charge, players get tokens from the 
pot.  In the second part of the game, the discharge, players get 
tokens from each other.  A player wins when he gets no token during 
the charge (many players may thus win), or when he first gets rid of 
his tokens during the discharge.

The charge or discharge is a sequence of sets.  In every set, each 
player at his turn plays a round against the dice, while the others 
wait.  The active player casts the dice up to three times; after every 
cast, he can put aside any number of dice, thus accumulating a 
combination.  Next players must cast dice as many times as the first 
player.\footnote{The order of players in the set matters, but I could 
not find definite rules for its determination.} End-of-round 
accumulated combinations, obtained by all players in the set, are 
ranked in the hierarchic order
\begin{multline}
  421\succ 111\succ 611\succ 666\succ 511\succ 555\succ 411\succ 
  444\succ 311\\
  \succ 333\succ 211\succ 222\succ 654\succ 543\succ 432\succ 321\succ
  665\succ\ldots 221,
  \label{ordhie}
\end{multline}
where $\succ$ means `higher than'.  The combinations, implicit in 
(\ref{ordhie}), are ordered as the numbers formed by their faces in a 
decreasing sequence: e.  g. $655\succ 654$.  The dominant combination 
421 and the dominated combination 221, known as ``n\'enette'', differ 
only by one die.  $fff$ is the $f$-brelan, $f11$ is the $f$-pair 
($f\neq 1$), $654,543,432,321$ are the sequences.

At end of set, the last\footnote{The adjective `last' is my own 
suggestion for automatic tie-breaking.} player who has got the lowest 
combination gets the number of tokens determined by table 
\ref{transf}, e.  g. if the highest combination is 411, then the last 
player with the lowest combination (whatever it is) gets 4 tokens.  
During the charge, tokens are taken from the pot, if possible.  When 
the pot is empty, the discharge begins, and tokens are now taken from 
the player who has got the highest combination.

\begin{table}
    \caption{token transfer function}
\label{transf}
\begin{center}
\begin{tabular}{|c|c|}
\hline highest combination&token number\\
\hline 421 &$10$\\
\hline 111 &$7$\\
\hline $f11,fff,f\neq 1$ &$f$\\
\hline sequence&$2$\\
\hline other&$1$\\
\hline 
\end{tabular}
\end{center}
\end{table}

\section{A Galton-Watson process in the 421 round\label{GW}}
 
Taking the genealogic point of view, each die is considered as an 
individual, dying after being cast, either without a child, in case of 
accumulation, or with a single child (itself indeed).  The child 
number being lower than one, the number of live dice $D_{j}$ (section 
\ref{dual}) decreases in time.  Moreover, the population becomes 
extinct after $J$ casts (or sooner).

A Galton-Watson process \cite{Harris} is obtained when the offspring 
of each individual is independent of others'.  With an optimal 
$\vec{d}^{*}$-goal strategy, the dice dying without children have 
their faces in $\vec{d}^{*}$, but the converse is not true.  For 
example, with $\vec{d}^{*}=221,\vec{d}_{1}=211,J>1$, the two dice 11 
have correlated offspring: one has a child if and only if the other 
has none.  Dice have independent offspring if and only if 
$\vec{d}^{*}$ is a brelan and the player is first.

I apply the Galton-Watson theory \cite[§6.2]{Parzen} to obtain the 
probability law of $D_{j}$, for an optimal $\vec{d}^{*}$-goal 
strategy, where $\vec{d}^{*}=D\vec{e}_{F}$ is the $F$-brelan.  Dice 
are indexed by $d=1\ldots D_{j}$ .  Let $Z_{d}\in\{0,1\}$ the number 
of children of the die indexed by $d$.
\begin{equation}
    D_{j}=\sum_{d=1}^{D_{j-1}}Z_{d}.
    \label{Dj++}
\end{equation}
The $Z_{d}$ are random variables, with the same law 
$q_{i}\equiv\mathcal{P}(Z_{d}=i)$, of generating function
\begin{equation}
    g(z)\equiv\<z^{Z_{d}}\>=q_{0}+q_{1}z,q_{0}=\frac{1}{F},q_{1}=1-q_{0}.  
    \nonumber
\end{equation}

The $Z_{d}$ are always independent if and only if $\vec{d}^{*}$ is a 
brelan and the player is first.  When this is true, from (\ref{Dj++}), 
the generating function of $D_{j}$, conditioned by $D_{j-1}$, is
\begin{equation}
    \<z^{D_{j}}|D_{j-1}=d\>=g(z)^{d}.
    \nonumber
\end{equation}
The generating function of $D_{j}$ is thus determined by
\begin{equation}
    g_{0}(z)=z^{D},g_{J}(z)=1,
    \nonumber
\end{equation}
\begin{multline}
    \forall(j,1\le j<J),g_{j}(z)\equiv\<z^{D_{j}}\>
    =\sum_{d=0}^{D}\<z^{D_{j}}|D_{j-1}=d\>\mathcal{P}(D_{j-1}=d)
    \\=\sum_{d=0}^{D}\mathcal{P}(D_{j-1}=d)g(z)^{d}
    =g_{j-1}\circ g(z).
    \nonumber
\end{multline}
By induction,
\begin{equation}
    g_{j}=g_{0}\circ g^{\circ j}.
    \nonumber
\end{equation}
The composition powers of the affine function $g$ are
\begin{equation}
    g^{\circ j}(z)=1-q_{1}^{j}+q_{1}^{j}z.
    \nonumber
\end{equation}
Therefore
\begin{eqnarray}
    g_{j}(z)&=&(1-q_{1}^{j}+q_{1}^{j}z)^{D},
    \nonumber\\
    \mathcal{P}(D_{j}=d)&=&\begin{pmatrix}D\\d\end{pmatrix}
    (1-q_{1}^{j})^{D-d}q_{1}^{jd}.
    \label{binome}
\end{eqnarray}

$D_{j}$ follows a binomial law, directly obtained by considering that 
a die dies when accumulated, or stays alive, with the probability 
$q_{1}$ per time step, independently of others: a Bernoulli process is 
obtained, with the law (\ref{binome}).  The interest of considering a 
Galton-Watson process is in the analogy with branching processes 
\cite{Harris}.

\section{Realization with \mma\label{mma}}

The present article is supported by \cite{421mma}, an open source 
software and data base in the \mma\ language \cite{Wolfram}, which, 
like \texttt{LISP}, is interpreted and allows functional and recursive 
treatments on arbitrary expressions, equivalent to trees.  The \mma\ 
frontend allows literate programming \cite{Knuth2} in the form of 
notebooks, gathering live code, outputs and comments, within a tree 
structure, that can be unfolded at will.

Combination manipulation differs slightly from list manipulation 
(since order does not matter in combinations) or ensemble manipulation 
(since repetitions are allowed in combinations).  A tool box is 
developed.  The numerical parameters $(D,F,J)$ are arbitrary, which 
realizes a scalable model, invaluable for development.  Fate trees are 
created recursively.  All fates converging to the same state at the 
same time are merged by indexing, so that the size grows only linearly 
with the depth $J$ and remains easily manageable for 
$(D,F,J)=(3,6,3)$.  In exchange, the computing time is increased and 
the history is lost, which allows to treat only time and state 
dependent utilities (as required in the 421 set).

Starting from the leaves of the fate tree, where utility is grafted, 
optimal strategies and expected utilities are build recursively, 
according to the $\max-\moy$ algorithm.  A utility-strategy tree is 
finally obtained, from which the strategy can be extracted, then piped 
into the $\str-\moy$ algorithm, a variation on $\max-\moy$, solving 
the Kolmogorov equation.

$\max-\moy$ produces the expected optimal one-goal strategies, 
Bernoulli or ratchet, depending on $D<F$, and dilemmas.  The result 
probabilities are computed, saved, and many properties are checked 
systematically.  Some result probabilities are checked by Monte Carlo 
simulations, with success.  The charts \ref{diaPro}, \ref{nonDiaPro}, 
\ref{nonDiaPro1}, \ref{nonDiaPro2}, \ref{nonDiaPro3} are generated 
automatically.  There is very little room for errors, and if there are 
any, they are traceable.

The goal identification heuristic policies are realized.  Their wrong 
decisions are pointed out.  They are exactly evaluated with 
$\str-\moy$, which is very slow, since it requires the computation of 
every heuristic decision in the fate tree, according to an algorithm 
actually longer and slower, for one decision, than the simple 
maximization in $\max-\moy$.  Obtaining the truth about heuristic 
policies is a lengthy task.

\begin{center}
    \textbf{Probability charts player's guide}
\end{center}    

p1, p2 mean first or next players.  In every box of a diagonal 
probability chart stands a column of the probabilities, ordered from 
top to bottom by growing delay, to obtain the goal written at head of 
line.  

In every box of a non-diagonal probability chart, stand two columns: 
at left, from top to bottom, the goal and the result; at right, the 
probabilities, ordered from top to bottom by growing delay, to obtain 
the result, with the goal in mind (and taking optimal decisions as 
determined by the ratchet).  Moreover, for easy access, the couples 
(goal, result) are represented in a square array, where heads of lines 
and columns are the respective representatives of goal and result, 
modulo face permutations (section \ref{1goalprob}).  The three-die 
representative 3X3 array is spread onto the three charts 
\ref{nonDiaPro1}, \ref{nonDiaPro2}, \ref{nonDiaPro3}, one for each 
goal class.

Here is an example for using non-diagonal charts.  Let the goal be 
$641$ and the result $652$.  The representatives of $641$ and $655$ 
are, separately and respectively, $123$ and $112$.  (Representatives 
are chosen so as to minimize the sum of their faces.)  The 
representative of the \emph{couple} $(641,655)$ is $(123,144)$.  $123$ 
takes us to chart \ref{nonDiaPro3} (the third line of the 
representative square array), whence $112$ takes us to the second 
column, $(123,144)$ to the third row, where finally are the 
probabilities to obtain, with the goal 641, the result 652, after one, 
two or three casts.


\begin{table}[tbp]
    \centering
    \caption{diagonal result probabilities}
    \vspace{1em}
    \includegraphics{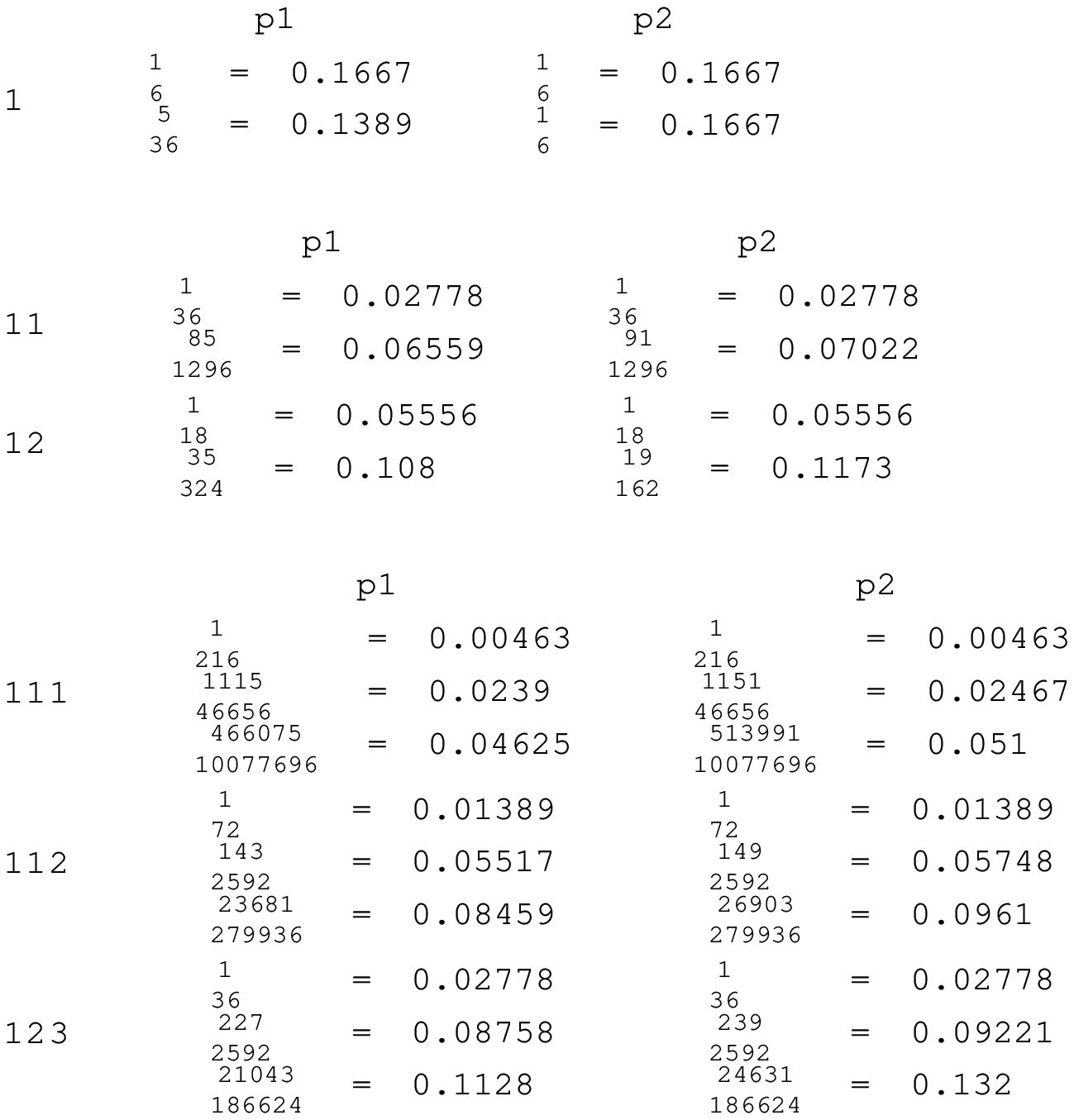}
    \label{diaPro}
\end{table}

\begin{table}[tbp]
    \centering
    \caption{first player's non-diagonal result probabilities}
    \vspace{1em}
    \includegraphics{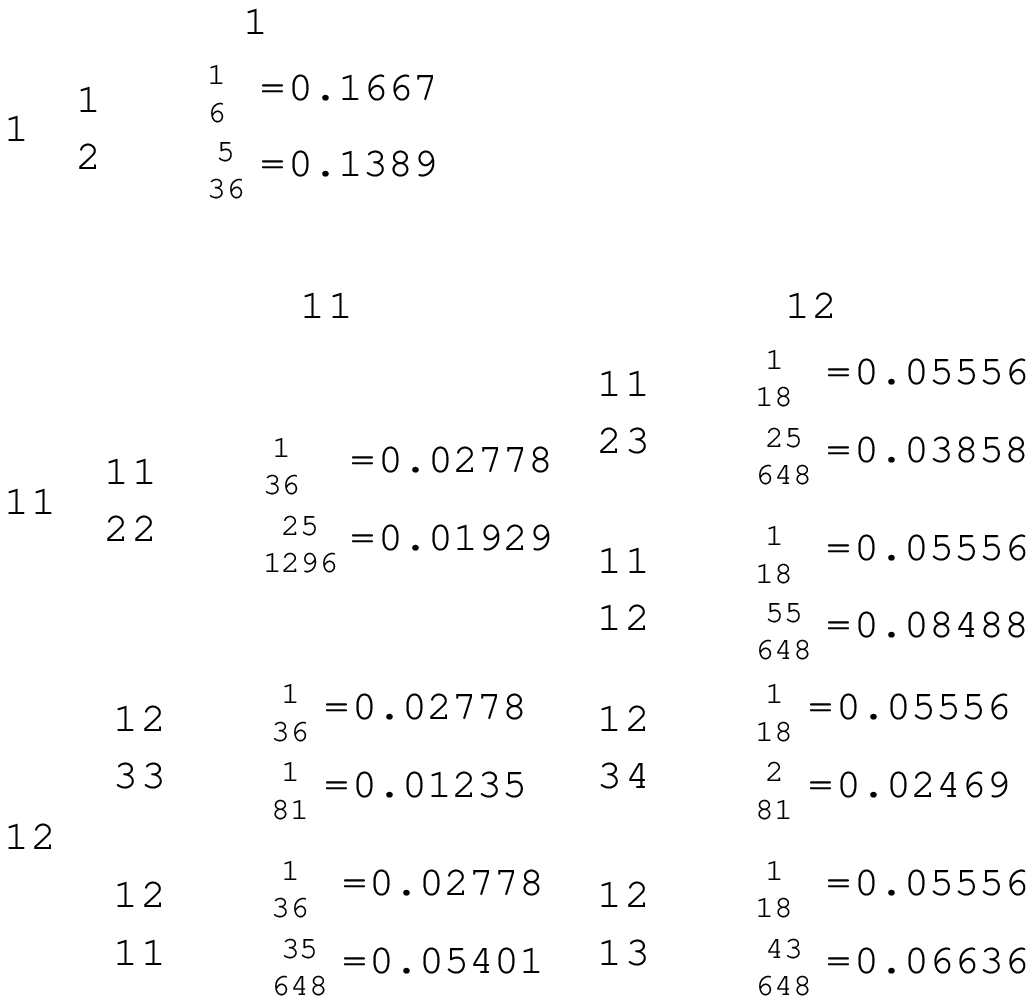}
    \label{nonDiaPro}
\end{table}

\begin{table}[tbp]
    \centering
    \caption{first player's non-diagonal result probabilities (1)}
    \vspace{1em}
    \includegraphics[scale=.95,angle=90]{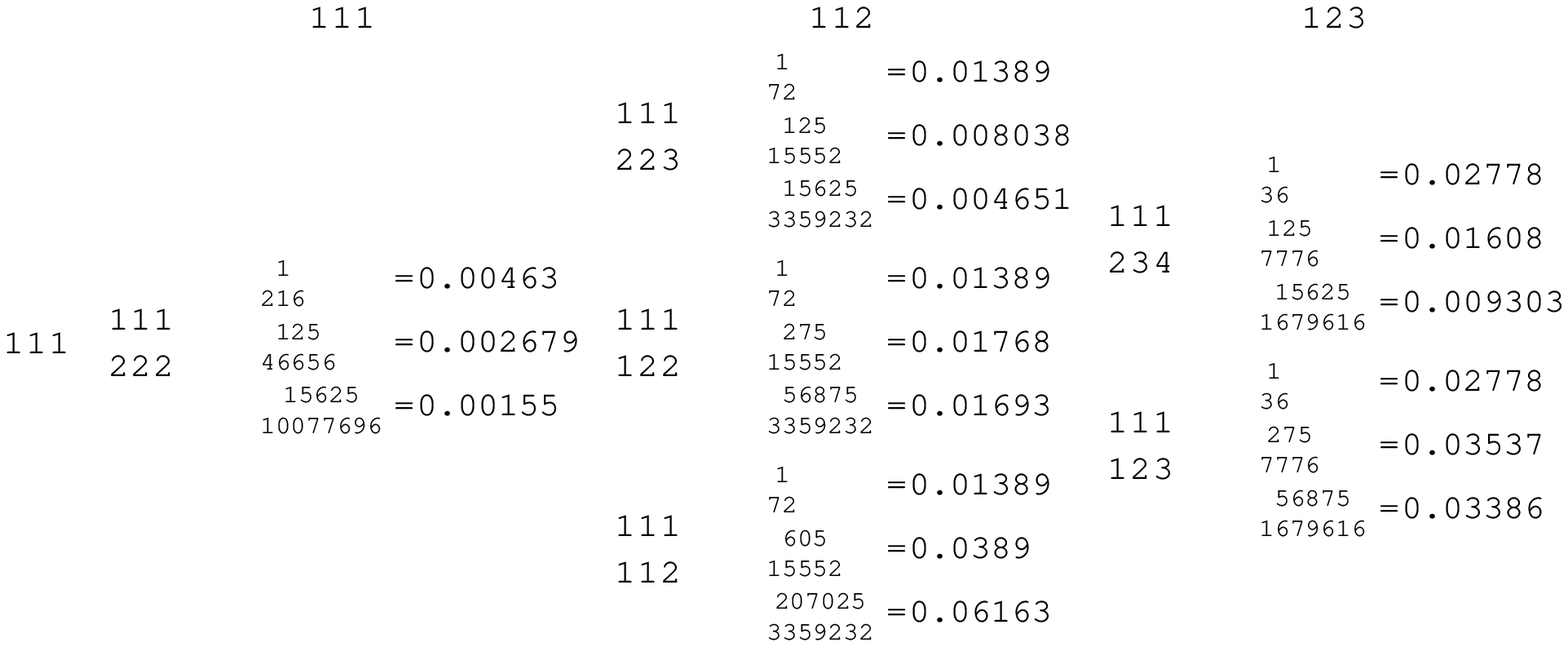}
    \label{nonDiaPro1}
\end{table}

\begin{table}[tbp]
    \centering
    \caption{first player's non-diagonal result probabilities (2)}
    \vspace{1em}
    \includegraphics[scale=.95,angle=90]{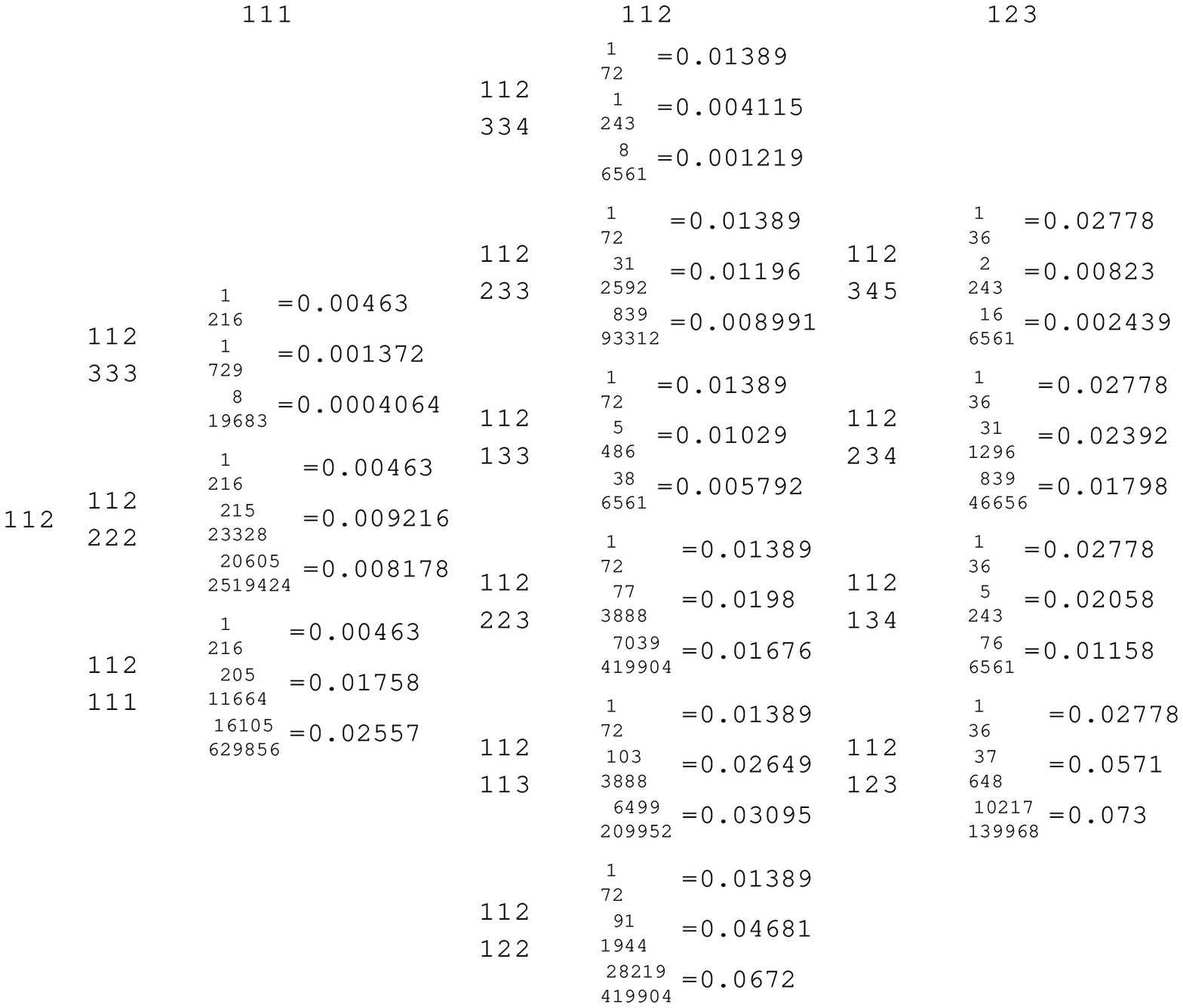}
    \label{nonDiaPro2}
\end{table}

\begin{table}[tbp]
    \centering
    \caption{first player's non-diagonal result probabilities (3)}
    \vspace{1em}
    \includegraphics[scale=.95,angle=90]{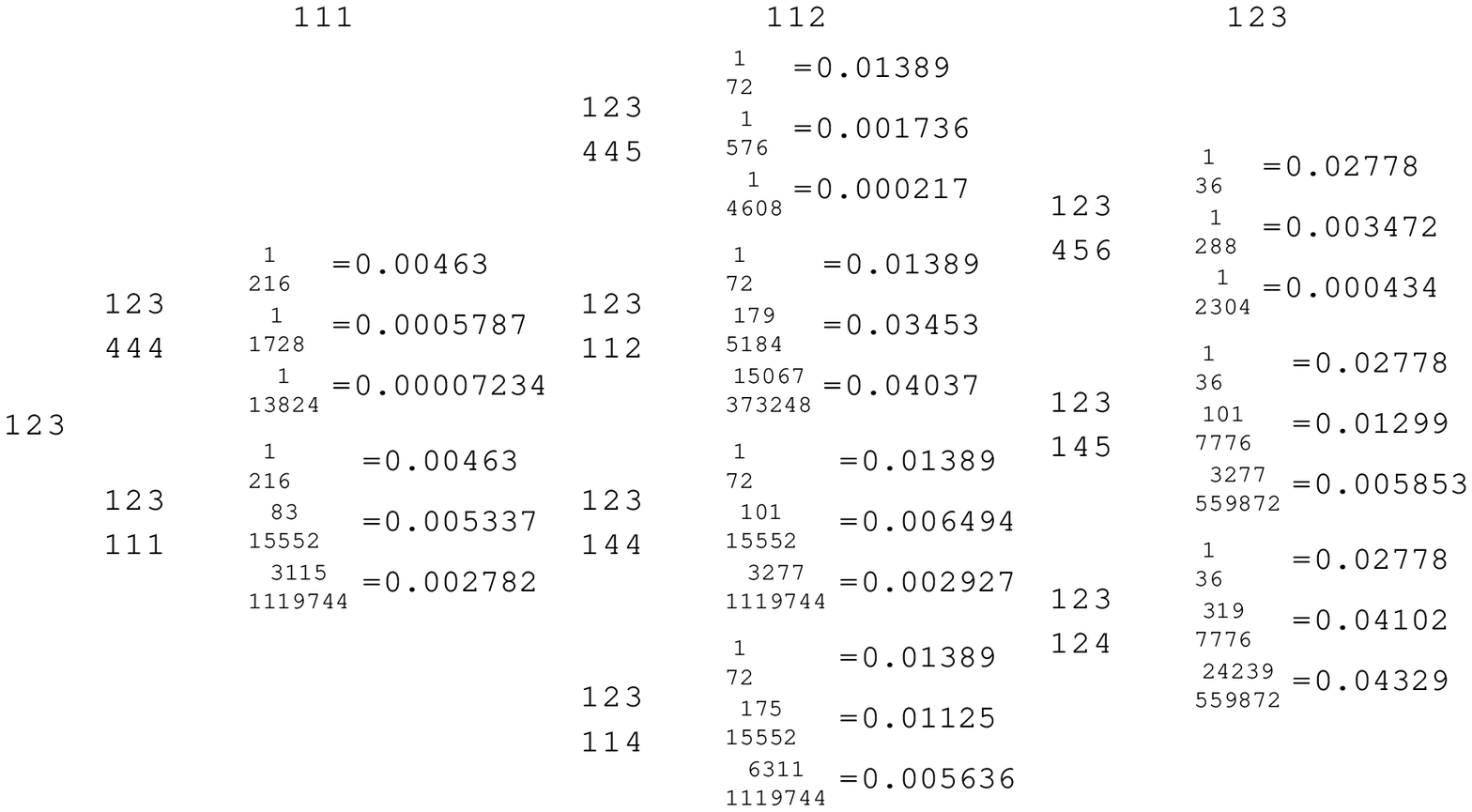}
    \label{nonDiaPro3}
\end{table}
\clearpage
\bibliographystyle{unsrt}
\bibliography{biblio}

\end{document}